\documentstyle[amscd,amssymb,verbatim,11pt]{amsart}

\theoremstyle{plain}
\newtheorem{Prop}{Proposition}[section]
\newtheorem{Thm}[Prop]{Theorem}
\newtheorem{Cor}[Prop]{Corollary}

\newtheorem{Lem}[Prop]{Lemma}

\theoremstyle{definition}
\newtheorem{Def}[Prop]{Definition}

\theoremstyle{remark}
\newtheorem{Rem}[Prop]{Remark}
\newtheorem{ex}[Prop]{Example}

\def\int{\mathop{\roman{int}}}

\def\1{^{-1}}

\def\tt{{\xi}}
\def\LL{{\mathcal L}}

\def\NN{{\mathbb N}}
\def\RR{{\mathbb R}}
\def\TT{{\mathcal S}}
\def\UU{{\mathcal U}}

\def\TL{\bar{\mathcal S}}

\errorcontextlines=0 \numberwithin{equation}{section}

\begin{document}
\title[Semigroup-controlled asymptotic dimension]{Semigroup-controlled asymptotic
dimension}

\author{J.~Higes}
\address{Departamento de Geometr\'{\i}a y Topolog\'{\i}a,
Facultad de CC.Matem\'aticas.
Universidad Complutense de Madrid.
Madrid, 28040 Spain}
\email{josemhiges@@yahoo.es}
\date{August 30, 2006}



\thanks{The author was supported by Grant AP2004-2494 from the Ministerio de Educacion y Ciencia, Spain.
He thanks  N. Brodskiy, J. Dydak, A. Mitra and J.M. Rodriguez
Sanjurjo for their support and helpful comments.}

\begin{abstract}
We introduce the idea of semigroup-controlled asymptotic
dimension. This notion generalizes the asymptotic dimension and
the asymptotic Assouad-Nagata dimension in the large scale. There
are also semigroup controlled dimensions for the small
scale and the global scale. Many basic properties of the
asymptotic dimension theory are satisfied by a
semigroup-controlled asymptotic dimension. We study how these new
dimensions could help in the understanding of coarse embeddings
and uniform embeddings. In particular we have introduced
uncountable many invariants under quasi-isometries and uncountable
many bi-Lipschitz invariants. Hurewicz type theorems are
generalized and some applications to geometric group theory are
shown.
\end{abstract}

\maketitle

\tableofcontents
\section{Introduction and preliminaries}

The aysmptotic dimension was introduced by Gromov in \cite{Gro
asym invar} as one important coarse invariant in the study of
geometric group theory. The idea behind its definition is to
analyze a metric space as a large scale object. An analogous
concept to the asymptotic dimension but for small scale objects
would be the uniform dimension introduced by Isbell in
\cite{Isbell book}. These two dimensions suggest the idea of
seeing a metric space as a global object and then we would get a
global definition of dimension. In \cite{Brod Dydak Higes Mitra2}(see also \cite{Protasov})
it was given such definition and it was analyzed some properties
of zero dimensional spaces.
\par The following notions associated to a cover $\UU =\{U_s\}_{s \in S}$ of a metric space $(X, d)$
are standard concepts in the
theory of metric spaces. They are used in many equivalent definitions of asymptotic dimension.\par
Let $\UU =\{U_s\}_{s \in S}$ be a cover of a metric space $(X, d_X)$, not necessarily open. Associated to this cover there
is a family a natural family of functions $\{f_s\}_{s \in S}$ with $f_s: X \to \RR_+$ defined by:

$f_s(x) := d_X(x, X \setminus U_s)$.

With such functions we can define:
\begin{itemize}
\item[-]{\itshape Local Lebesgue number} $\LL_x(\UU)$ of $\UU$ at $x \in X$: $\LL_x(\UU) := sup\{f_s(x)|S \in S\}$.
\item[-]{\itshape Global Lebesgue number} $\LL(\UU):= inf\{\LL_x(\UU)| x \in X\}$.
\item[-]{\itshape Local $s$-multiplicity} $s-m_{\UU}(x)$ of $\UU$ at $x \in X$ is
defined as the number of elements of $\UU$
that intersect
$B(x, s)$.
\item[-]{\itshape Global $s$-multiplicity} $s-m(\UU):= sup\{s-m_{\UU}(x) | x \in X\}$. If $s =0$ then the $0$-multiplicity
will be called
{\itshape multiplicity} of $\UU$ and it will be noted by $m(\UU)$.
\end{itemize}

Given a family of subsets $\UU$ of a metric space $(X,d_X)$ it is
said that $\UU$ is {\itshape $C$-bounded} with $C >0$ if $diam(U)
\le C$ for every $U \in \UU$. If $\UU$ is $C$-bounded for some
$C>0$ it is said that $\UU$ is {\itshape uniformly bounded} and if
$d_X(U, V) > s$ for every $U, V \in \UU$,  it is said also that
$\UU$ is {\itshape $s$-disjoint} with $s > 0$. \par

A definition of asymptotic dimension could be the following:
\begin{Def}
We will say that a metric space $(X, d)$ has asymptotic dimension
at most $n$ (notation: $asdim X \le n$) if there is an $s_0$ such
that for every $s\ge s_0$ there exists an uniformly cover(colored
cover) ${\mathcal{U}} = \bigcup_{i =1}^{n+1} {\mathcal{U}}_i$ so
that each ${\mathcal{U}}_i$ is $s$-disjoint.
\end{Def}

From this definition we can deduce that there exist a function $f: \RR_+ \to \RR_+$ with $\lim_{x\to \infty} f(x) = \infty$
such that each $s$-disjoint $n+1$ colored cover is $f(s)$-bounded. If we restricted the range of functions $f$ allowed to the
linear ones(or asymptotically linear ones) we would get the notion of asymptotic Assouad-Nagata dimension, also called
linear controlled asymptotic dimension or asymptotic dimension with Higson property. Many results have
appeared recently associated to this dimension, see \cite{Brod-Dydak-Higes-Mitra}, \cite{Brod Dydak Higes Mitra2},
\cite{Dran Smith}, \cite{Dran-Zar} and
\cite{Hur}, for example. Here is the definition:

\begin{Def}\label{Nagata Definition}
We will say that a metric space $(X, d)$ has asymptotic
Assouad-Nagata dimension at most $n$ (notation: $asdim X \le n$)
if there is a $s_0 >0 $ and a $C>0$ such that for every $s\ge s_0$
there exist a $Cs$-bounded cover(colored cover) ${\mathcal{U}} =
\bigcup_{i =1}^{n+1} {\mathcal{U}}_i$ so that each
${\mathcal{U}}_i$ is $s$-disjoint.
\end{Def}

If we just change the condition of $s \ge s_0$ by $s \le s_0$ we will get the notion of capacity dimension introduced
by Buyalo
 in \cite{Buyalo1}. The asymptotic Assouad-Nagata dimension and the capacity dimension are just the
 large scale and the small scale
versions of the Assouad-Nagata dimension introduced by Assouad in
\cite{Assouad}. The Assouad-Nagata dimension is a bi-Lipschitz
invariant. Lang and Schlichenmaier proved in \cite{Lang-Sch Nagata
dim} that the Assouad-Nagata dimension is in fact a quasisymmetry
invariant. Many other interesting properties of the Assouad-Nagata
dimension appeared in that work. \par

The relations among the capacity dimension, the asymptotic Assouad-Nagata dimension and
the Assouad-Nagata dimension were showed
in \cite{Brod-Dydak-Higes-Mitra}.

In this paper we generalize all these notions using the concept of
semigroup-controlled asymptotic dimension, semigroup-controlled
small scale dimension and semigroup-controlled global dimension.
The idea is to modify the range of functions allowed to control
the size of the colored covers in definition \ref{Nagata
Definition}. These functions will have a semigroup structure. In
section \ref{Control semigroups} we will study some basic
properties of this kind of semigroups. In the next section we will
define the notion of semigroup-controlled dimension and we will
show that many basic properties of asymptotic dimension (see
\cite{Bell Dran Beld}) are satisfied for these new dimensions.
Section \ref{Large scale small scale} is dedicated to relate the
large scale, small scale and global dimensions following the same
ideas of \cite{Brod-Dydak-Higes-Mitra}. In section \ref{Non
equivalent semigroup-controlled dimensions} we will prove that
this generalization is not trivial i.e. we have introduced
uncountable many different dimensions. The last two sections are
dedicated to important topics. In section \ref{Invariance of
dimension} we will study the types of maps between metric spaces
for which these new dimensions are invariant. In particular we
will prove that all these dimensions are quasi-isometric
invariants in the large scale theory or bi-Lipschitz invariants in
the global case. Last section generalizes the results of
\cite{Hur} about the Hurewicz type theorem. Some applications to
geometric group theory are obtained as corollaries.

\section{Control semigroups}\label{Control semigroups}

In this work we will consider properties ${\mathcal{P}}(s)$
depending on positive real numbers $s \in \RR_+$. We will say that
{\itshape ${\mathcal{P}}(s)$ is satisfied in a neighborhood of
$\infty$ (respectively in a neighborhood of $0$) if there is a
$s_0$ such that ${\mathcal{P}}(s)$ is true for every $s \ge
s_0$(resp. $s \le s_0$).}

Properties that are satisfied in a neighborhood of $\infty$ will
be called large scale properties or asymptotic properties,
properties that are satisfied in a neighborhood of $0$ will be
called small scale properties. If a property ${\mathcal{P}}(s)$ is
satisfied in a neighborhood of $\infty$ and in a neighborhood of
$0$ we will say that ${\mathcal{P}}(s)$ is a global property. \par
In this work we will usually give proofs and statements for the
large scale case. The small scale case and the global case can be
usually done  using a dual reasoning and in many cases they will
be left to the reader.

Next definition describes the type of functions that are going to control the dimension.

\begin{Def}\label{Def Control Functions}
Let $f: \RR_+ \to \RR_+$ be an increasing continuous function
\begin{itemize}
\item[a.] We will say that it is a large scale(or asymptotic)
dim-control function if $f(x) \ge x$ in a neighborhood of $\infty$
and $f(\infty) = \infty$ i.e. $lim_{x \to \infty} f(x) = \infty$.
\item[b.] If we require that $f(x) \ge x$ in a neighborhood of $0$
and $f(0) =0$ we will say that $f$ is a small scale dim-control
function. \item[c.]If a function is a small scale dim-control
function and a large scale dim-control function then we will say
that it is a global dim-control function.
\end{itemize}
\end{Def}
If a large scale dim-control function $f$ is equal to some linear
function in a neighborhood of $\infty$ we will say that it is a
 {\itshape large scale(or asymptotically) linear dim-control function.}
Analogously it can be defined the notion of small scale linear dim-control function and
global linear dim-control function.

Now we define the notion of large scale(or asymptotic) control
semigroup. This notion is the main concept of this paper.

\begin{Def}
Let $\TT$ be a set of asymptotic dim-control functions. We will say that
$\TT$ is a large scale(or asymptotic) control semigroup if the following properties are
satisfied:
\begin{enumerate}
\item Every asymptotically linear dim-control function is in
$\TT$(linear condition). \item For every pair of functions $g_1,
g_2 \in \TT$ we have $g_1 \circ g_2 \in \TT$(semigroup condition).
\end{enumerate}

We define the notion of small scale and global control semigroup for a set of small scale(resp. global)
dim-control function
analogously. We will note it by $\tt$ (resp. $\bar{\TT}$)
\end{Def}

\begin{Rem}
In this work we will refer to asymptotic control semigroup just as
control semigroups unless it were necessary to remark the large
scale condition.
\end{Rem}

Here are the main examples of control semigroups.
\begin{ex}
The set of all asymptotic dim-control functions is a control semigroup. It will be called
the uniform control semigroup and we will note it by ${\itshape U}$ (${\itshape u}$ and $\bar{\itshape U}$ for
the small scale and the global case respectively)

The set of all asymptotically linear dim-control functions is a control semigroup.
It will be called the Nagata(or the linear) control semigroup and it will be noted by $\mathcal{N}$
(${\itshape n}$ and $\bar{\mathcal{N}}$ for
the small scale and large scale case respectively).
\end{ex}

\begin{Rem}
Let $K$ be a set of asymptotic dim-control functions. From the theory of semigroups we know that
such set generates a control semigroup $L(K)$ i.e.
the intersection of all control semigroups $\TT$ so that
$K \subset {\TT}$. It will be
called the {\itshape control semigroup generated by $K$}. Note that
${\mathcal{N}} = L(\{x\})$. We have also obtained many examples of control semigroups
using this procedure.
\end{Rem}

In the set of all control semigroups we can define a partial order.
\begin{Def}
Given two control semigroups $\TT_1$ and $\TT_2$, we will say that
$\TT_2$ is  finer than $\TT_1$ (notation: ${\TT_1} \preceq \TT_2$)
if for every asymptotic dim-control function $f \in \TT_2$ there
is a dim-control function $g \in \TT_1$ such that $f \le g$ in a
neighborhood of $\infty$. If $\TT_1 \preceq \TT_2$ and $\TT_2
\preceq \TT_1$ we will say that both control semigroups are
equivalent and we will note it by $\TT_2 \approx \TT_1$.
\end{Def}

The following result is easy to check and it is left as an exercise.

\begin{Prop} Let $\TT$ be a control semigroup then:
\begin{enumerate}
\item ${ \itshape U } \preceq {\TT}\preceq {\mathcal{N}} $.
\item $ L(K \cup {\TT}) \preceq
\TT$ with $K$ is any set of asymptotic dim-control functions.
\end{enumerate}
\end{Prop}

\begin{Rem}
As a consequence of the second statement in the previous
proposition we obtain that $L({\TT_2}\cup {\TT_1}) \preceq
{\TT_2}$ and $L({\TT_2}\cup {\TT_1}) \preceq {\TT_1}$ for every
pair of control semigroups. Then the set of all control semigroups
with the relation $\preceq$ is a directed set.
\end{Rem}

Now we show an example of two different equivalent control semigroups. In the following section we could
see how this example allows us to give another definition of asymptotic Assouad-Nagata dimension:

\begin{ex}\label{Ex Piecewise And Nagata}
Let $\{C_i\}_{i=1}^n$ be a finite set of constants with $C_i \ge
1$. For each finite set of those constants we can define an
asymptotic dim-control function $f_{\{C_i\}}$ as a continuous
piecewise linear function built with linear functions of slope
$C_i$ such that $f_{\{C_i\}}(\infty) = \infty$. It can be easily
checked that the set of all functions of the form $f_{\{C_i\}}$ is
a control semigroup. We will note it by $\mathcal{PL}$. Clearly we
have ${\mathcal{PL}} \preceq \mathcal{N}$. Take a $f_{\{C_i\}} \in
{\mathcal{PL}}$. Let $C$ be the maximum of the $C_i$ then the linear function $f(x) = Cx$
satisfies $f(x) \ge f_{\{C_i\}}(x)$ for all $x$ in a neighborhood of $\infty$.
So we have proved
that ${\mathcal{N}} \preceq {\mathcal{L}}$ and both semigroups are
equivalent.
\end{ex}

The following is an easy consequence of the semigroup theory and it will be used later.

\begin{Prop}\label{generation of dimfamilies}
Let $K$ be a set of dim-control functions.  If $g$ is a
dim-control function $g$ of $L(K)$ then there exist a finite
sequence of dim-control functions $\{f_i\}_{i=1}^n$ such that:

$g = f_1 \circ f_2 \circ ...\circ f_n$\\
where each $f_i$ belongs to ${\mathcal{N}} \cup K$.
\end{Prop}

\section{Semigroup-controlled dimensions: basic properties}\label{Semigroup-controlled dimensions: basic properties}

Now, using the notion of control semigroup we can give the
definition of semigroup-controlled asymptotic dimension. It
generalizes the notions of Assouad-Nagata asymptotic dimension and
asymptotic dimension.

\begin{Def}
Let $\TT$ be a control semigroup. We will say that a metric space
$(X, d)$ has $\TT$-controlled asymptotic dimension at most $n$
(notation: $asdim_{\TT} X \le n$) if there is an $f \in \TT$ such
that for every $s$ in some neighborhood of $\infty$ there exist a
cover(colored cover) ${\mathcal{U}} = \bigcup_{i =1}^{n+1}
{\mathcal{U}}_i$ so that each ${\mathcal{U}}_i$ is $s$-disjoint
and $f(s)$-bounded. \par

A metric space is said to have $\TT$-controlled asymptotic dimension $n$ if it has $\TT$-controlled asymptotic
dimension at most $n$
and for every $k < n$
it does not happen that $asdim_{\TT} X \le k$.
\end{Def}

\begin{Rem} \begin{itemize}
\item[a.] The function $f \in \TT$ for which $asdim_{\TT} X \le n$
will be called {\itshape $(n, \TT)$-dimensional control function
of $X$}. Such notion will become very important in the last
section. \item[b.] The small scale $\tt$-controlled dimension will
be noted by $smdim_{\tt} X$. \item[c.] Special remark is needed
for the definition of global dimension. The strictly analogous
definition to the Assouad-Nagata dimension would be the following:

\par $(X, d)$ has $\TL$-controlled global dimension at most $n$
(notation: $dim_{\TL} X \le n$) if there is a $f \in \TL$
such that for every $s$ there is a cover(colored cover)
${\mathcal{U}} = \bigcup_{i =1}^{n+1} {\mathcal{U}}_i$ so that each
${\mathcal{U}}_i$ is $s$-disjoint and $f(s)$-bounded.\\

\end{itemize}
\end{Rem}

We will say that two control semigroups $\TT_2$ and $\TT_1$
are {\itshape $dim$-equivalent} (notation: $asdim_{\TT_2} \equiv
asdim_{\TT_1}$) if for every metric space $(X,d)$ we have
$asdim_{\TT_2} X = asdim_{\TT_1} X$.

Next proposition and corollaries justify our definition of $\preceq$.

\begin{Prop} \label{families finer minor dim}
Let $\TT_1$ and $\TT_2$ be two control semigroups. If
${\TT_1}\preceq {\TT_2}$ then for every metric space $X$
we have $asdim_{\TT_1} X \le asdim_{\TT_2} X$.
\end{Prop}

\begin{pf}
Given $f \in \TT_2$ take $g \in \TT_1$ such that $g(s) \ge f(s)$ in a neighborhood of $\infty$,
then any $f(s)$-bounded family $\UU$
of subsets of $X$ is $g(s)$-bounded with $s$ in a neighborhood of $\infty$.
\end{pf}

\begin{Cor}
Let $\TT$ be a control semigroup. For every metric space $X$ we have
$asdim_{\itshape U} X \le asdim_{\TT} X \le asdim_{\mathcal{N}} X$.
\end{Cor}

\begin{Cor}\label{Equivalent Means Dim-Equivalent}
If two control semigroups are equivalent then they are $dim$-equivalent, i.e. they define the same asymptotic dimension.
\end{Cor}

\begin{ex}
Note that $dim_{\itshape \bar{U}} X$ and $asdim_{\itshape U} X$ are the uniform dimension defined in
\cite{Brod Dydak Higes Mitra2}
and the asymptotic dimension
of \cite{Gro asym invar} respectively.
We have also that $asdim_{{\mathcal N}} X$, $smdim_{\itshape n} X$ and
$dim_{\bar{\mathcal N}} X$ are the asymptotic dimension with Higson property(see for example \cite{Roe lectures}),
the capacity dimension(see \cite{Buyalo1})
and the Assouad Nagata dimension(introduced in \cite{Assouad}
) respectively.
\end{ex}

Note that applying \ref{Equivalent Means Dim-Equivalent} to the example
\ref{Ex Piecewise And Nagata} we get another way of defining the Assouad-Nagata asymptotic dimension.

The following proposition covers some basic properties of a
semigroup-controlled asymptotic dimension. The proofs are highly
similar to that ones for Assouad Nagata dimension and asymptotic
Assouad Nagata dimension. We quote between brackets the works where the
analogous proofs can be founded.

\begin{Prop}\label{Prop Bus}
Let $(X, d_X)$ and $(Y, d_Y)$ be metric spaces and $\TT$ a control semigroup. Then it is satisfied:
\begin{enumerate}
\item $asdim_{\TT} A \le asdim_{\TT} X$ for every $A \subset X$. (This is trivial)
\item $asdim_{\TT} X \times Y \le asdim_{\TT} X + asdim_{\TT} Y$ (\cite{Hur} or \cite{Lang-Sch Nagata dim})
\item If $X = A \cup B$ then $asdim_{\TT} X = max \{asdim_{\TT} A, asdim_{\TT} B\}$
(\cite{Hur} or \cite{Lang-Sch Nagata dim}).
\item The following conditions are equivalent: (\cite{Bell Dran Beld})
\begin{itemize}
\item[a.] $asdim_{\TT} X \le n$. \item[b.] There is a $f \in \TT$
such that for every  $s$ in some neighborhood of $\infty$ there
exist a cover ${\mathcal{U}}$ with $s-m(\UU)\le n+1$ and
$f(s)$-bounded. \item[c.] There is a  $g \in \TT$ such that for
every $s$ in some neighborhood of $\infty$ there exist a cover
${\mathcal{U}}$ with $m(\UU) \le n+1$, ${\mathcal{L}}({\UU}) \ge
s$ and $g(s)$-bounded \item[d.] There is a  $h \in \TT$ such that
for every $\epsilon$ in some neighborhood of $0$ there is a map
$\epsilon$-Lipschitz $p: X \to K^n$ with $K^n$ a $n$-dimensional
simplicial complex  such that the family $p^{-1}(st_v)$ is
$h(1/\epsilon)$-bounded.
\end{itemize}
\end{enumerate}
\end{Prop}

\begin{Rem}
Note that in the proof of the third property it is used the semigroup condition.
For the second and fourth properties it is necessary
that given $g_1, g_2 \in \TT$ there exists a $g_3 \in \TT$ such that $g_1 + g_2 \le g_3$ in a neighborhood of $\infty$.
Using the
semigroup condition and the linear condition define
$g_3 = 2\cdot g_1\circ g_2$ if $g_1(x) > x$ and $g_2(x) > x$ in a neighborhood of $\infty$.
\end{Rem}

\section{Large scale and small scale dimensions}\label{Large scale small scale}
The aim of this section is to study how the large scale, small scale and the global dimensions are related. Many of these
results are based on \cite{Brod-Dydak-Higes-Mitra}.

Given a global control semigroup we can see it just as a large scale(or small scale) control semigroup. This is the idea behind
next definition.

\begin{Def}
Let $\TL$ be a global control semigroup we will define the large
scale (resp. small scale) truncated semigroup of $\TL$ as the
semigroup of all functions $g$ for which there exist a dim-control
function $f \in {\TL}$ with $g(x) = f(x)$ in a neighborhood of
$\infty$ (resp. in a neighborhood of $0$). We will note it by:
$Trunc^{**}({\TL})$ (resp. $Trunc_{**}({\TL})$).
\end{Def}

Now we present some kind of inverse operation of truncation. Given a small scale control semigroup $\tt$ and a large scale
control semigroup $\TT$ we want to create a global control semigroup $\TL$ associated to those ones.

\begin{Def}
We define the linked set of $\tt$ and $\TT$ as the set of all
continuous increasing functions $g$ for which there exist a small
scale dim-control function $g_1$ with $g_1 \in {\tt}$ and a large
scale dim-control function $g_2$ with $g_2 \in \TT$ such that
$g(x) = g_1(x)$ in a neighborhood of $0$ and $g(x) = g_2(x)$ in a
neighborhood of $\infty$. It will be noted by $Link({\tt}, \TT)$.
\end{Def}

Clearly $Trunc^{**}({\TL})$ and $Trunc_{**}({\TL})$ are small scale and large scale control semigroups respectively. Next
proposition shows we have the same property for linked sets.

\begin{Prop}\label{Pasted Families Are Dim}
Let $Link({\tt}, \bar{\TT})$ be a linked set then it is a global control semigroup.
\end{Prop}
\begin{pf}
The linear condition is trivial. Now let $f, g$ be two dim-control
functions of $Link({\tt}, \TT)$ and let $(f_1, f_2)$, $(g_1, g_2)$
its small scale and large scale associated functions. we have that
$f(g(x)) = f_1(g_1(x))$ in a neighborhood of $0$ and $f(g(x)) =
f_2(g_2(x))$ in a neighborhood of $\infty$ then as $f_1
\circ g_1 \in \tt$ and $f_2 \circ g_2 \in \TT$ the semigroup
condition is satisfied.
\end{pf}

The relationship between truncation and linking is given in the following result:

\begin{Prop} \label{Truncate And Paste Is Stupid}
Let $\TL$ be a global control semigroup then:\par
$Link(Trunc_{**}({\TL}), Trunc^{**}({\TL})) \approx \TL$. \\
Conversely $Trunc^{**}(Link(\tt, \TT))\approx \TT$ and $Trunc_{**}(Link(\tt, \TT))\approx \tt$.
\end{Prop}
\begin{pf}
Let ${\mathcal{P}}$ be the global control semigroup $Link(Trunc_{**}({\TL}), Trunc^{**}({\TL}))$.
Clearly by ${\TL}\subset {\mathcal{P}}$ we have ${\mathcal{P}}\preceq {\TL}$.

Now let $g$ be a dim-control function in $\mathcal{P}$.
There exist two dim-control functions
$g_1, g_2 \in \TL$ such that $g(x)= g_1(x)$ in a neighborhood of $0$
and $g(x) = g_2(x)$ in a neighborhood of $\infty$. Let $x_1 \le x_2$ be two positive numbers such that
$g(x) = g_1(x)$ if $x \le x_1$ and $g(x) = g_2(x)$ if $x\ge x_2$. Let $M$ be the maximum of $g$ in $[x_1, x_2]$
and let $m$ be the minimum of $g_1 \circ g_2$ in the same interval. Suppose $m < M$ then there exist
a $C >1$ such that $Cm \ge M$.
Define the function $g'(x) = C(g_1(g_2(x)))$. We have $g'(x) \ge g(x)$ for every $x$ and then
${\TL}\preceq {\mathcal P}$.

The converse is obvious.
\end{pf}

Next definition is the key to connect large scale, small scale and global dimensions, see \cite{Brod-Dydak-Higes-Mitra}.

\begin{Def} Let $\TL$ be a global control semigroup and let $X$ be any metric space.
We will say that $(X,d)$ has $\TL$-microscopic controlled dimension at most $n$
and we will note by $m-dim_{\TL} X \le n$
if the metric space $(X, d'_1= min(d, 1))$ has
$\TL$-controlled dimension at most $n$.
In a similar way we will say that a metric space
$(X,d)$ has $\TL$-macroscopic controlled dimension at most $n$
($M-dim_{\TT} X \le n$) if $dim_{\TL} (X, d''_1) \le n$ with $d''= max(1, d)$.
\end{Def}

\begin{Lem}\label{C Does Not Matter}
Let $(X,d)$ be a metric space and let $\TL$ be a global control semigroup.
If for every $c >0$ we define the metrics
$d'_c =min (c, d)$ and $d''_c = max(c,d)$
then
$m-dim_{\TL} X = dim_{\TL} (X, d'_c)$ and $M-dim_{\TL} X = dim_{\TL} (X, d''_c)$.
\end{Lem}

The proof of this result is given in \cite{Brod-Dydak-Higes-Mitra}. In such proof the authors used that the Assouad-Nagata
dimension is invariant under Lipschitz functions.
We give another proof without using this fact.

\begin{pf}
We just do the proof for the microscopic case. The macroscopic case is similar.
Suppose without loss of generality that $c < 1$.
Let $s$ be any positive number. If $s \ge c$ then pick the cover ${\mathcal{U}} = X$.
Assume that $s < c$ then there is a cover ${\mathcal{U}}$ with $s$-Lebesgue number in $(X, d'_1)$,
multiplicity at most $m-dim_{\TL} X +1$ and $f(s)$-bounded
for some dim-control function $f \in \TL$. But as $d'_1(x,y) = d'_c(x,y) = d(x,y)$ if $d(x,y) \le c$
then such cover satisfy $\LL(\UU) \ge s$ in $(X, d'_c)$ and it is $f(s)$-bounded.
We have proved $dim_{\TL} (X, d'_c) \le m-dim_{\TT} X$.
The remaining case is similar.
\end{pf}

\begin{Cor}
Let $\TL_1$ be a global control semigroup. If $X$ is a bounded metric
space then $m-dim_{\TL_1} X = dim_{\TL_1} X$. If $X$ is
a discrete metric space then $M-dim_{\TL_1} X =
dim_{\TL_1}X$.
\end{Cor}

Next lemma shows how the microscopic dimension of a global control
semigroup is greater or equal than the semigroup-controlled
asymptotic dimension associated to the large scale truncated
semigroup.

\begin{Lem}\label{Macroscopic Dimension And Covering}
Let $X$ be a metric space and let $\TL$ be a global control semigroup then the
following properties are equivalent:
\begin{enumerate}
\item $M-dim_{\TL} X \le n$
\item There is a function $f \in \TL$ such that for all $s$ in a neighborhood of $\infty$ there is a colored cover
${\mathcal{U}} = \bigcup_{i =1}^{n+1} {\mathcal{U}}_i$ with each
${\mathcal{U}}_i$ $s$-disjoint and $f(s)$-bounded.
\end{enumerate}
\end{Lem}

The proof is almost equal to the proof of Lemma $2.7.$ of \cite{Brod-Dydak-Higes-Mitra}. It will be left to the reader.

Using the same reasoning we can get the analogous Lemma for the microscopic case:

\begin{Lem}
Let $X$ be a metric space and $\TL$ a global control semigroup then the
following properties are equivalent:
\begin{enumerate}
\item $m-dim_{\TL} X \le n$
\item There is a function $f \in \TL$ such that for all $s$ in a neighborhood of $0$
 there is a colored cover ${\mathcal{U}} = \bigcup_{i =1}^{n+1} {\mathcal{U}}_i$ with each ${\mathcal{U}}_i$ $s$-disjoint and $f(s)$-bounded.
\end{enumerate}
\end{Lem}

Next lemma could be considered some kind of converse of the previous ones.
\begin{Lem}
Let $\TL$ be a global control semigroup. For every metric space $X$ we
have $dim_{\TL}X \le n$ if and only if $m-dim_{\TL}
X \le n$ and $M-dim_{\TL} X \le n$.
\end{Lem}
\begin{pf}
The necessary condition is obvious by lemma ~\ref{Macroscopic
Dimension And Covering}. Let us prove the sufficient condition.
Suppose $m-dim_{\TL} X \le n$ and $M-dim_{\TL}X \le n$. Without
loss of generality we can assume that the $f$ associated to the
bounds of the microscopic covers and the macroscopic covers is the
same, if not take the composition. Let $s$ be a positive real
number. We want to find a dim-control function $g \in \TL$ and a
colored covering $\mathcal{U}$ $s$-disjoint and $g(s)$-bounded. It
is clear that if $s >1$ or $f(s) < 1$ the result is obvious.
Assume that $s\le 1$ and $f(s) \ge 1$. Pick $s_0 = f^{-1}(1)$ and
define the function $g(x) = f(2f(x))$. So if $s_0 \le s \le 1$
take a colored covering $\mathcal{U}$ of $(X, d''_1)$ so that it
is $2f(s)$-disjoint and $f(2(f(s)))$-bounded.
\end{pf}

Combining all the results of this section we get that the global dimension can be obtained just studying
the dimension in a neighborhood of $0$ and in a neighborhood of $\infty$.

\begin{Thm}\label{Truncate And Dim Is The Same}
Let $\TL$ be a global control semigroup. For every metric space $X$ we
have that $smdim_{Trunc_{**}({\TL})} X = m-dim_{\TL} X$ and
$asdim_{Trunc^{**}({\TL})} X = M-dim_{\TT} X$.

In the other hand given $\tt$ a small scale control semigroup and $\TT$ a large scale control semigroup then:
\begin{center}
 $dim_{Link({\tt}, {\bar{\TT}})} X = max \{smdim_{\tt} X,
asdim_{\TT} X\}$.\\
\end{center}

\end{Thm}
\begin{pf}
By the two previous lemmas we get $M-dim_{\TL} X \ge
asdim_{Trunc^{**}({\TL})} X$. Suppose $asdim_{Trunc^{**}({\TL})} X
\le n$. Then for every $s \ge s_0$ there is a cover
${\mathcal{U}}^s = \bigcup_{i =1}^{n+1} {\mathcal{U}}_i^s$ that is
$s$-disjoint and $f(s)$-bounded with $f \in Trunc^{**}({\TL})$.
That means that there is a $f_1 \in \TL$ such that $f(x) = f_1(x)$
if $x \ge x_1$ for some $x_1$. Let $s'_0$ be the maximum of $x_1$
and $s_0$. We have that for every $s \ge s'_0$ there is a cover
${\mathcal{U}}^s = \bigcup_{i =1}^{n+1} {\mathcal{U}}_i^s$
$s$-disjoint and $f_1(s)$-bounded. Applying ~\ref{Macroscopic
Dimension And Covering} we get the result. The microscopic case is
analogous.

For the second statement just note:

$dim_{Link({\tt}, {\TT})} X =
max\{m-dim_{Link({\tt}, {\TT})} X, M-dim_{Link({\tt}, {\TT})} X\}=$\\
\par
$ max\{smdim_{Trunc_{**}(Link({\tt}, \TT))}X, asdim_{Trunc^{**}(Link({\tt}, {\TT}))} X\}$\\\\
So the result follows from the second statement of \ref{Truncate
And Paste Is Stupid} and \ref{families finer minor dim}.
\end{pf}

\section{Non equivalent semigroup-controlled dimensions}\label{Non equivalent semigroup-controlled dimensions}

Let $\overline{ASDIM}$(respectively $\overline{SMDIM}$) be the quotient set of all large scale (resp. small scale)
control semigroups with the equivalence relation $\equiv$ defined in section
\ref{Semigroup-controlled dimensions: basic properties}. In this section we will estimate the cardinality of
$\overline{ASDIM}$ and $\overline{SMDIM}$.

Let $\Omega$ be the set of all the countable ordinals($\Omega_0$) union the first uncountable ordinal. This set is
uncountable and it has a natural well order. We will prove that there exist two order preserving maps
$i_L: \Omega \to {\overline{ASDIM}}$ and $i_s: \Omega \to {\overline{SMDIM}}$. As usual
we will prove the results for the large
scale case. The small scale case will be left to the reader.

As a first step we will estimate the cardinality of the set of
all large scale(resp. small scale) control semigroups
modulo the equivalence relation $\approx$. They will be noted by $ASDIM$ and $SMDIM$.

The next lemmas are necessary to show that there exist at least
countable many semigroups.

\begin{Lem}\label{function is bigger}
Let $g:[a_i, a_{i+1}] \to [g(a_i), g(a_{i+1})]$ be an increasing
continuous function defined in some interval and let $f(a_i),
f(a_{i+1})$  be any pair of points that satisfies $f(a_i) \ge
g(a_i), f(a_{i+1}) \ge g(a_{i+1})$ then there is a continuous
increasing  function $\bar{f}$ defined in the same interval such
that $\bar{f} \ge g$ with $\bar{f}(a_i) = f(a_i)$ and
$\bar{f}(a_{i+1}) = f(a_{i+1})$ .
\end{Lem}
\begin{pf}
Define the function $f_{aux}(x) = g(x) + f(a_i)-g(a_i)$. This is an
increasing continuous function in the interval. If
$f_{aux}(a_{i+1})=g(a_{i+1})$ then $f_aux = \bar{f}$. Otherwise
let $z$ be the greatest point in $[a_i, a_{i+1}]$ such that the
segment $\overline{f(a_i)f(a_{i+1})}$ intersects the graph of
$g$. Such point exists by continuity. Then define:

$\bar{f}(x) = \begin{cases} f_{aux}(x) \text{ if } x\in [a_i, z]\\
\frac{f(a_{i+1})-f_{aux}(z)}{a_{i+1}-z}(x-z)+f_{aux}(z) \text{  if
}
x\in (z, a_{i+1}]\end{cases}$\\
It is clear that this function satisfies the requirements of the
lemma.
\end{pf}

\begin{Lem}\label{function is superbigger}
Let $f$ be a large scale dim-control function such that it is strictly bigger than
the identity in a neighborhood of $\infty$. Then there exist a large scale
dim-control function $g$ such that for every $n \in \mathbb{N}$  $g(x) > f^n(x)$ in some neighborhood of $\infty$.
\end{Lem}

\begin{pf}
As the function is strictly bigger than the identity we have that
there is an increasing sequence $x_n \to \infty$ such that if $x
\ge x_n$ then $f^n(x)
> f^i(x)$ for all $i\le n-1$. Using ~\ref{function is bigger}
define in each interval $[x_n, x_{n+1}]$ a function $g_n$ so that
$g_n(x_n) = f^n(x_n)$, $g_n(x_{n+1}) = f^{n+1}(x_{n+1})$ and
$g_n(x) \ge f^n(x)$. Paste all these functions and we get $g$ the
required function.
\end{pf}

Using a similar reasoning we can give the following lemma:

\begin{Lem}\label{function is superbigger2}
Let $f$ be a small scale dim-control function such that it is strictly bigger than
the identity in a neighborhood of $0$. Then there exist a
small scale dim-control function $g$ such that for every $n \in \mathbb{N}$  $g(x) > f^n(x)$ in some neighborhood of $0$.
\end{Lem}

Next lemma shows that there are at least countable many non
equivalent large scale control semigroups.

\begin{Lem}\label{infinite families}
There exist a sequence $\{\TT_i\}_{i\in \mathbb{N}}$ of large scale
control semigroups with $\TT_i \prec \TT_{i-1}$ and
$\TT_i = {L}(\{f_i\} \cup {\TT}_{i-1})$ where $f_i$
is an asymptotic dim-control function such that for every $g \in {\TT}_{i-1}$
$f_i(x) > g(x)$ in a neighborhood of $\infty$.
\end{Lem}
\begin{pf}
Take $\TT_1 = {\mathcal N}$. Let $f:
\mathbb{R}_+ \to \mathbb{R}_+$ be the dim-control function defined by:

$f(x)= \begin{cases} x^2 \text{ if } x \in [1, \infty)\\ x \text{
otherwise }\end{cases}$.

It is clear that for every asymptotically linear dim-control function $g$ there is a point
$x_0 \in \mathbb{R}_+$ such that $f(x) > g(x)$ if $x \ge x_0$ then
we have ${\TT_2}=L(\{f\} \cup \mathcal{N}) \prec {\mathcal{N}}
= \TT_1$.

Suppose we have constructed a sequence of control semigroups with
$\TT_n \prec...\prec \TT_1$ so that for each
control semigroup $\TT_i$ there is an asymptotic dim-control function $f_i$ that
${\TT_i }= L(\{f_i\} \cup \TT_{i-1})$ and for every
dim-control function $g$ of $\TT_{i-1}$ there exists a $x_0 \in
\mathbb{R}_+$ such that $f_i(x)
> g(x)$ if $x > x_0$. Now apply lemma ~\ref{function is superbigger}
to $f_n$ in order to get a dim-control function $f_{n+1}$ so that for
every $j\in \mathbb{N}$, $f_{n+1}(x) > f_n^j(x)$ if  $x>x_0$ for
some $x_0$. We claim that $f_{n+1}$ satisfies the same property
for all $g \in \TT_n$. Let $g$ be a dim-control function in
${\TT}_n$. Using ~\ref{generation of dimfamilies} we have
that:

$g = f_1 \circ f_2 \circ ...\circ f_p$\\
For some functions $f_i$ in $\{f_n\} \cup \TT_{n-1}$. For every
$f_i$ there exist an $x_i$ and a $j_i$ such that $f_n^{j_i}(x) >
f_{i}(x)$ if $x \ge x_{i}$. Let $x'_0$ be the maximum of all
$x_{i}$ then we have that $g(x) < f^{{\sum_{i=1}^p j_i}}_n(x)$. We
have obtained that $g(x) < f_n^j(x)$ for some $j$ if $x \ge x_0$
and by the method we have built $f_{n+1}$ we have $f_{n+1}(x)
> g(x)$ in a neighborhood of $\infty$.

\end{pf}

Note that for getting $\TT_2$ we just need an asymptotic dim-control function $f$ that were strictly greater than any
asymptotically linear dim-control function in a neighborhood of $\infty$.

Doing a dual reasoning we can get:

\begin{Lem}\label{infinite families2}
There exist a sequence $\{{\tt}_i\}_{i\in \mathbb{N}}$ small scale
control semigroups with $\tt_i \prec \tt_{i-1}$ and
${\tt}_i = L(\{f_i\} \cup {\tt}_{i-1})$ where $f_i$
is a dim-control function such that for every $g \in {\tt}_{i-1}$ $f_i(x) > g(x)$ in a neighborhood of $0$.
\end{Lem}

We have proved that the sets $SMDIM$ and $ASDIM$ are at least countable.
The lemmas of above suggest the following definition:

\begin{Def}
We will say that a large scale(small scale) control semigroup $\TT$ (resp. $\tt$) is mono-bounded at  $\infty$
(respectively at $0$)
if there exist an asymptotic dim-control function $f$ such that for every  $g \in \TT$ (resp $g \in \tt$) we have
 $f(x) > g(x)$ in a neighborhood of
$\infty$ (resp. in a neighborhood of $0$).
The function $f$ will be called the bound function of $\TT$ (of $\tt$) at $\infty$ (resp. at $0$).
\end{Def}

Using the sequence generated in \ref{infinite families}
we can build the control semigroup generated by such sequence. Such semigroup will be mono-bounded
and then we can begin again a similar process as in \ref{infinite families}.
This is the idea of the next two lemmas.

\begin{Lem}\label{Everything Is Mono}
Let $\{\TT_i\}_{i \in \mathbb{N}}$ be a sequence of
large scale control semigroups mono-bounded at $\infty$ such that  $\TT_i \prec \TT_{i-1}$ and
$\TT_i = L(\{f_i\} \cup {\TT}_{i-1})$ where $f_i$ is a bound function of $\TT_{i-1}$
then the large scale control semigroup given by $L(\bigcup_{i = 1}^{\infty} \TL_i)$
is mono-bounded at $\infty$.
\end{Lem}

\begin{pf}
The reasoning is similar to the previous ones. \par
Firstly we note that $\TT:= L(\bigcup_{i = 1}^{\infty} \TT_i) = \bigcup_{i = 1}^{\infty} \TT_i$.

We will prove that $\TT$ is mono-bounded at $\infty$.
We have that there is a sequence of points
 $\{a_i\}_{i\in \mathbb{N}}$ with $a_i \to \infty$  and $a_{i+1}> a_i$ such that  $f_{i+1}(x) > f_i(x)$
for all $x \ge a_{i+1}$. Using ~\ref{function is bigger} we can build a function $f$
such that $f(a_i) = 2f_i(a_i)$ and
$f(x) \ge f_i(x)$ in $[a_i, a_{i+1}]$. Let $g$ be a dim-control function of ${\TT}_i$ by hypothesis there is a
$x_0$
such that
$f_{i+1}(x) > g(x)$ for every $x \ge x_0$ and as we have $f(x) > f_{i+1}(x)$ for every $x \ge a_{i+2}$ then
 $\TT$ is mono-bounded at $\infty$ and $f$ is a bound function.
\end{pf}

Here is the lemma for the small scale:
\begin{Lem}\label{Everything Is Mono2}
Let $\{\tt_i\}_{i \in \mathbb{N}}$ be a sequence of
small scale control semigroups mono-bounded at $0$ such that  $\tt_i \prec \tt_{i-1}$ and
$\tt_i = L(\{f_i\} \cup {\tt}_{i-1})$ where $f_i$ is the bound function of $\tt_{i-1}$
then the small scale control semigroup given by $L(\bigcup_{i = 1}^{\infty} \tt_i)$
is mono-bounded at $0$.
\end{Lem}

Now as $\bigcup_{i=1}^{\infty} \TT_i$ is mono-bounded we can apply
again \ref{infinite families} and generate another strictly
decreasing sequence of control semigroups and such semigroups will
be mono-bounded. Applying repeatedly the arguments of
\ref{Everything Is Mono} and \ref{infinite families} (or
\ref{Everything Is Mono2} and \ref{infinite families2} in the
small scale case) we get:

\begin{Cor}\label{Embeddings Of Omega}
There are two injective maps
$i_L: \Omega \to ASDIM$, $i_s: \Omega \to SMDIM$ such that
$i_L(\beta) \prec i_L(\alpha)$ ($i_s(\beta) \prec i_s(\alpha)$) for every pair of ordinals with $\alpha < \beta$.
\end{Cor}
\begin{pf}
By the previous reasoning it is obvious that there exist an
injective map
 $i_L^0: \Omega_0 \to {ASDIM}\setminus \{{\itshape U}\}$ that it is order preserving.
Moreover we have ${\itshape U} \prec i_L^0(\alpha)$ for every
$\alpha \in \Omega_0$. Define:\par
$i_L(\alpha) = \begin{cases} i_L^0(\alpha) \text{ if } \alpha \in \Omega_0 \\
{\itshape U} \text{ Otherwise }.\end{cases}$\\
The small scale case is analogous.
\end{pf}


Now we will prove that the control semigroups given by the map of \ref{Embeddings Of Omega} are not
$dim$-equivalent.

\begin{Lem}\label{Infinite dimensions}
Let $\TT$ be a mono-bounded control semigroup and let $f$ be a bound function at $\infty$
of $\TT$. Then there exist a metric space $X$
such that $0 = asdim_{L(\{f\} \cup \TT)} X < asdim_{\TT} X$.
\end{Lem}

\begin{pf}

Clearly $f$ is strictly greater than the identity. Let
$\{a_i\}_{i\in \NN} \subset \RR_+$ be the sequence defined by $a_0 = 0$ and $a_i
= f(n_i)$ with $n_0 = 0$ and $n_i$ being the smallest natural number such that $n_i > n_{i-1}$ and
$a_{i-1} + n_i < f(n_i)$ for every $i \in {\mathbb{N}}$. Define the sequence of sets $\{C_i\}_{i=1}^{\infty}$ as:
\begin{center}
$C_i := \{a_{i-1}\}\cup \{a_i\}\cup \{\{a_{i-1}+m\cdot n_i\}_{m\in \NN} \cap [a_{i-1}, a_i]\}$.\\
\end{center}
Take the set $X := \bigcup_{i= 1}^{\infty} C_i$ we claim that this set satisfies the required properties. \par
Note that $U_0^{n_k} := \bigcup_{i=1}^k C_i$ is $n_k$-connected and $diam(U_0^{n_k}) = f(n_k)$. Also note that if
$s \in [n_k, n_{k+1})$
then the $s$-connected components of $X \setminus U_0^{n_k}$ have cardinality
at most two and $d(X\setminus U_0^{n_k}, U_0^{n_k}) > n_{k+1}$. Hence $asdim_{L(\{f\} \cup \TT)} X \le 0$.\par
Now if $asdim_{\TT} X \le 0$ with $g$ as $(0, \TT)$-dimensional control function take $x_0$ such that $f(x) > g(x)$ if
$x \ge x_0$. Then for every $n_k \ge x_0$ we have that every $n_k$-disjoint cover contains the $n_k$-connected set
$U_0^{n_k}$ but as
$diam(U_0^{n_k}) = f(n_k) > g(n_k)$ we get a contradiction.
\end{pf}

For the small scale the lemma is similar.
The proof will follow the steps of the previous lemma but we will construct $X$ as
a convergent sequence to $0$ nor to $\infty$. The details are left to the reader.

\begin{Lem}\label{Infinite dimensions2}
Let $\tt$ be a mono-bounded small scale control semigroup and let $f$ be a bound function at $0$
of $\tt$. Then there exist a metric space $X$
such that $0 = smdim_{L(\{f\} \cup \TT)} X < smdim_{\TT} X$.
\end{Lem}

As a trivial
corollary of the previous results we get the following theorem that gives us an estimation of the cardinalities of
$\overline{ASDIM}$ and $\overline{SMDIM}$.

\begin{Thm}\label{First Theorem}
There exist two injective maps: \par
$i_L: \Omega \to \overline{ASDIM}$\par
$i_s: \Omega \to \overline{SMDIM}$\\
such that each $i_L(\alpha)$ (resp. $i_s(\alpha)$) is  mono-bounded
at $\infty$(resp. at $0$) if $\alpha \in \Omega_0$. Moreover we have that
$i_L(\beta) \prec i_L(\alpha)$ $(i_s(\beta) \prec i_s(\alpha))$ for every pair of ordinals such that $\alpha < \beta$.
\end{Thm}

\section{Maps between metric spaces and dimension}\label{Invariance of dimension}
The aim of this section is to study when a
semigroup-controlled dimension of a metric space is invariant
under modifications of the metric. The type of modifications
allowed is described in the next definition.
\begin{Def}
Consider a set $\mathcal{Q}$ of positive real continuous
increasing functions $\{\rho_{\lambda}\}$.
\begin{enumerate}
\item If $\rho_{\lambda}(\infty) = \infty$ for every $\rho_{\lambda} \in \mathcal{Q}$
we will say that $\mathcal{Q}$ is a large scale metric perturbation set or just metric
perturbation set.
\item If
$\rho_{\lambda}(0) = 0$ for every $\rho_{\lambda} \in \mathcal{Q}$ we will say that it is a small scale metric perturbation set.
\item If $\mathcal{Q}$ is a large scale and a small scale metric perturbation set we will say that
it is a global perturbation
set.
\end{enumerate}
\end{Def}

Using this definition we can define the notion of coarse embedding in the following way:

\begin{Def} \label{Definition Of Embedding}
Let $\mathcal{Q}$ be a metric perturbation set.
Let $X, Y$ be two metric spaces and
let $f: X \to Y$ be a function. We will say that $f$ is a $\mathcal{Q}$-coarse embedding
if there are two functions $\rho_-$ and
$\rho_+$ with $\rho_-^{-1}, \rho_+ \in \mathcal{Q}$ such that:

$\rho_{-}(d_X(x, y)) \le d_Y(f(x), f(y)) \le \rho_{+}(d_X(x, y))$ for every pair of points $x, y \in X$.\\
In particular if a function $f$ satisfies the previous inequality we can say that it is a $(\rho_{-}, \rho_+)$-coarse embedding.
Functions $\rho_-$ and $\rho_+$ are called dilatation  and contraction respectively.
\end{Def}

Analogous definitions are given for the small scale and global
theory. Note that a small scale embedding and a global embedding
are always injective. but a coarse embedding is not necessarily
injective.

\begin{Rem}\label{Remmark To Cite}
\begin{enumerate}
\item Let ${\mathcal{Q}}$ be the set of all positive continuous
increasing functions $\rho_{\lambda}$ with $\rho_{\lambda}
(\infty) = \infty$ , a large scale $\mathcal{Q}$- coarse embedding
is called a coarse embedding.
\item  If ${\mathcal{Q}}$ is the set of all linear functions then
a ${\mathcal{Q}}$-coarse embedding is called a
quasi-isometric embedding. See ~\cite{Gro asym invar}. \item If
${\mathcal{Q}}$ is the set of all positive continuous increasing
functions $\rho_{\lambda}$ with $\rho_{\lambda} (0) = 0$
  then a  small scale $\mathcal{Q}$-embedding is called uniform embedding
~\cite{Isbell book}. Also it is known as a uniformly continuous
embedding.
\end{enumerate}
\end{Rem}

The following proposition gives a sufficient condition for a semigroups controlled dimension to be invariant
under a metric perturbation set.

\begin{Prop} \label{Function Strange Is Useful}
Let $\mathcal{Q}$ be a metric perturbation set.
Let $\TT_1, \TT_2$ be two control semigroups. Suppose that
given any dim-control function $f$ with $f \in \TT_2$ and given any two functions $\rho_-, \rho_+$ with
$\rho_-^{-1}, \rho_+ \in \mathcal{Q}$
there exist a dim-control function $g \in \TT_1$ that satisfies
$\rho_{-}^{-1}\circ f \circ \rho_+ \le g$ in a neighborhood of $\infty$.
Then if $F: X \to Y$ is a $\mathcal{Q}$-embedding between
metric spaces $X, Y$ we have that  $asdim_{\TT_1}X \le asdim_{\TT_2} Y$.
\end{Prop}
\begin{pf}
Suppose $asdim_{\TT_2} Y \le n$. Let $t >0$ be a positive number
in a neighborhood of $\infty$. Take a cover $\mathcal{U}$ in Y of
multiplicity at most $n+1$, with Lebesgue number at least
$\rho_+(t)$ and uniformly bounded by $f(\rho_+(t))$ with $f \in
{\TT_2}$. Let $\mathcal{V}$ be the covering given by
${\mathcal{V}} = F^{-1}({\mathcal{U}})$. Note that  the
multiplicity of ${\mathcal{V}}$ is $n+1$. Also we have that if $d_X(x, y)
\le t$ then $d_Y(F(x), F(y)) \le \rho_+(t)$ and the Lebesgue
number of $\mathcal{V}$ is at least $t$. Now if we have two points
$x, y$ that are in the same set $V \in {\mathcal{V}}$ we have:

$d(x,y) \le \rho_-^{-1}(d(F(x),F(y)) \le \rho_-^{-1}(f(\rho_+(t)))$\\
By hypothesis there is a $g \in {\TT_1}$ such that $\rho_{-}^{-1}\circ f \circ \rho_+ \le g$. Hence $\mathcal{V}$ is
$g(t)$-bounded and we have proved $asdim_{\TT_1}X \le asdim_{\TT_2} Y$.
\end{pf}
It is clear that the proposition of above also works for small scale and global embeddings.

Now given a control semigroup $\TT$ we would like to find a metric
perturbation set $\mathcal{Q}$ for which the associated
semigroup-controlled asymptotic dimension $asdim_{\TT}$ is
invariant under any $\mathcal{Q}$-embedding, i.e. if $f: X \to Y$
is a $\mathcal{Q}$-embedding then $asdim_{\TT}X \le asdim_{\TT}
Y$. This is the idea behind next definition.

\begin{Def}
Let $\TT$ be a control semigroup. We define the metric
perturbation set $\Sigma(\TT)$ generated by $\TT$ as the set of
all continuous increasing functions $\rho_{\lambda}$ with and
$lim_{x\to \infty} \rho_{\lambda}(x) = \infty$ for which there
exist a dim-control function $g\in \TT$ such that $g \ge
\rho_{\lambda}$ in a neighborhood of $\infty$
\end{Def}

\begin{Rem}
Note that the semigroup condition allow us to show that if
$f: X \to Y$ is a $\Sigma(\TT)$-coarse embedding and $g: Y\to Z$ is a
$\Sigma(\TT)$-coarse embedding then $g\circ f: X \to Z$ is a $\Sigma(\TT)$-coarse embedding.
\end{Rem}

Using \ref{Function Strange Is Useful} it is easy to check that $asdim_{\TT}$ is invariant
under a $\Sigma(\TT)$-embedding.

\begin{Cor}\label{Main Theorem}
Let $\TT$ be any control semigroup. Let $X, Y$ be any two metric spaces. If there exist a $\Sigma(\TT)$-coarse embedding
$f: X \to Y$ then $asdim_{\TT} X \le asdim_{\TT} Y$.
\end{Cor}
\begin{pf}
Just note that the conditions of \ref{Function Strange Is Useful} are satisfied.
\end{pf}

Clearly the analogous corollary is true for global scale and small scale dimensions.

\begin{Rem} Using ~\ref{Main Theorem} and the results of \ref{First
Theorem} we have found uncountable many invariants under
quasi-isometries or uncountable many invariants under bi-Lipschitz
equivalences for the global case. Note also that if a metric
perturbation set $\mathcal{Q}$ is mono-bounded by $f$ then using a
reasoning similar as in \ref{First Theorem} but taking $\TT_1 =
L(\{f\})$ in \ref{infinite families}  and applying again \ref{Main
Theorem} we get that there exist uncountable many
$\mathcal{Q}$-invariants.
\end{Rem}

For the zero dimensional case we have the following nice theorem. The proof appears in ~\cite{Brod Dydak Higes Mitra2}.

\begin{Thm}\label{Dim Zero Is Ultrametric}
Let $\TL$ be  a global control semigroup. For every metric space $(X,d_X)$ with
$dim_{\TL} X \le 0$ there exist an ultrametric  $d_u$ in $X$
such that the identity $id: (X,
d_X) \to (X, d_u)$ is a $\Sigma(\TL)$-embedding.
\end{Thm}

Following the ideas of ~\cite{Brod Dydak Higes Mitra2} we have
also that every ultrametric space can be embedded via a
bi-Lipschitz function in $L_{\omega}$. The construction of
$L_{\omega}$ is the following: \par Let $S$ be a countable set and
fix an element $s_0 \in S$. The set $L_{\omega}$ is the subset of
all sequences $\bar{x} = \{x_n\}_{n \in \mathbb{Z}}$ with $x_n \in
S$. We will say that a sequence is in $L_{\omega}$ if there exist
a $k$ such that  $x_n = s_0$ for every $n < k$. The distance
between two elements of $L_{\omega}$ is given by  $d(\bar{x},
\bar{y}) = 3^{-m}$ if $m \in \mathbb{Z}$ is the minimum index such
that $x_m \ne y_m$. We have:

\begin{Cor}
Let $\TL$ be a global control semigroup. For every metric space $(X,
d_X)$ with $dim_{\TT} X \le 0$ there exist a $\Sigma(\TT)$-embedding $f: X \to L_{\omega}$.
\end{Cor}

The idea of a $\Sigma(\TT)$-coarse embedding suggests us an analogous
definition of a $\Sigma(\TT)$-coarse map or global $\Sigma(\TL)$-map.
That definitions would be the generalization of coarse maps(or asymptotically Lipschitz maps) and uniform maps.
This notion
will be very important in the next section.

\begin{Def}
Let $\TT$ be a control semigroup. We will say that a function $f$ between metric spaces $f: X \to Y$ is a
$\Sigma(\TT)$-coarse map if there is a
$\rho \in \Sigma(\TT)$ such that:

$d_Y(f(x), f(y)) \le \rho(d(x,y))$ for every $x, y \in X$.
\end{Def}

\section{Hurewicz type theorems}\label{Hurewicz type theorems}
The aim of this section is to generalize the results of \cite{Hur}.  In order to get such generalization we need to
modify the notion of an $m$-dimensional control function of a function $f: X \to Y$ between metric spaces. That is
is the idea behind the next definition.

\begin{Def}
Let $\TT_1$ and $\TT_2$ be two control semigroups.
We will say that a function $D^{(\TT_1,\TT_2)}: \RR_+ \times \RR_+ \to \RR_+$ is a $(\TT_1, \TT_2)$-function if for every
$\rho_1 \in \TT_1$ and for every $\delta \in \TT_2$ then there is a function $\rho_2 \in \TT_1$ such that
$D^{(\TT_1, \TT_2)}(\rho_1(x), \delta(x)) \le \rho_2(x)$ in a neighborhood of $\infty$.
\end{Def}
The definitions for the small scale case and the global case are obtained as usual.

Now the following definition generalize the definition $4.4.$ of \cite{Hur}. We recall the following definition
of the cited work:  Given a function $f: X \to Y$ between
metric spaces a subset $A \subset X$ is said to be {\itshape $(r, R)$-bounded} if for any
$x, y \in A $ with $d_X(x, y) \le r$ then
$d_Y(f(x), f(y)) \le R$.

\begin{Def}
Let $\TT_1$, $\TT_2$ be any control semigroups. Given a function
$f: X \to Y$ between metric spaces and given $m \ge 0$, an
{\itshape $m$-dimensional $(\TT_1, \TT_2)$- control function} of
$f$ is  a $(\TT_1, \TT_2)$-function $D^{(\TT_1,\TT_2)}_f: \RR_+
\times \RR_+ \to \RR_+$ such that for all $r_X >0$ and $R_Y >0$ in
a neighborhood of $\infty$ any $(\infty, R_Y)$-bounded subset $A$
of $X$ can be expressed as the union of $m+1$ sets whose
$r_X$-components are $D^{(\TT_1,\TT_2)}_f(r_X, R_Y)$-bounded.
\end{Def}

The generalization of $4.6.$ of \cite{Hur} is similar:

\begin{Def}
Let $\TT_1$, $\TT_2$ be any control semigroups. Given a function
$f: X \to Y$ between metric spaces and given $k\ge m+1 \ge 1$, an
{\itshape $(m,k)$-dimensional $(\TT_1, \TT_2)$- control function}
of $f$ is  a $(\TT_1, \TT_2)$-function $D^{(\TT_1,\TT_2)}_f: \RR_+
\times \RR_+ \to \RR_+$ such that for all $r_X >0$ and $R_Y >0$ in
a neighborhood of $\infty$ any $(\infty, R_Y)$-bounded subset $A$
of $X$ can be expressed as the union of $k$ sets $\{A_i\}_{i=1}^k$
whose $r_X$-components are $D^{(\TT_1,\TT_2)}_f(r_X, R_Y)$-bounded
so that any $x \in A$ belongs to at least $k-m$ elements of
$\{A_i\}_{i=1}^k$.
\end{Def}
Definitions for the small scale and global case are obvious.

With these generalizations the analogous theorem (for the large scale case) of $4.9.$ in \cite{Hur} would be the following.

\begin{Prop}\label{Quasi Hur}
Let $\TT_1$, $\TT_2$ be any control semigroups and let $k = m+n+1$ where $m, n \ge 0$.
Suppose $f: X \to Y$ is a $\Sigma(\TT_2)$-coarse function of metric spaces and $asdim_{\TT_2} Y \le n$. If there is an
$(m,k)$-dimensional $(\TT_1, \TT_2)$-control function $D^{(\TT_1,\TT_2)}_f$ of $f$ then:
\begin{center}
$asdim_{\TT_1} X \le m + n$\\
\end{center}
\end{Prop}

The proof is almost the same as in \cite{Hur}. Note that in the proof the semigroup condition is used strongly.
For small scale case and the global case there are analogous theorems. We just give the global one:

\begin{Prop}
Let $\TL_1$, $\TL_2$ be any global
control semigroups and let $k = m+n+1$ where $m, n \ge 0$.
Suppose $f: X \to Y$ is a global $\Sigma(\TL_2)$-function of metric spaces and $dim_{\TL_2} Y \le n$. If there is an
$(m,k)$-dimensional $(\TL_1, \TL_2)$-control function $D^{(\TL_1,\TL_2)}_f$ of $f$ then:
\begin{center}
$dim_{\TT_1} X \le m + n$\\
\end{center}
\end{Prop}

Now note that as a trivial consequence of proposition $4.7.$ of \cite{Hur} we have:
\begin{Prop}
Let $f: X \to Y$ be a function of metric spaces and $m \ge 0$. Suppose that there exist an
$m$-dimensional $(\TT_1, \TT_2)$-
control function of $f$ with $\TT_1$, $\TT_2$ two control semigroups. Then there exist an
$(m,k)$-dimensional $(\TT_1, \TT_2)$-
control function of $f$ for every $k >m+1$
\end{Prop}

We can define the dimension of a function between metric spaces in the following way:

\begin{Def}
Given $\TT_1$, $\TT_2$ two control semigroups and given a function $f: X\to Y$ of metric spaces we define the
$(\TT_1, \TT_2)$-asymptotic dimension $asdim_{(\TT_1, \TT_2)}(f)$ of $f$ as the minimum of
$m$ for which there is an $m$-dimensional
$(\TT_1, \TT_2)$-control function.
\end{Def}

Finally we have as a corollary of \ref{Quasi Hur}

\begin{Thm}\label{Hur La}
If $f: X \to Y$ is an $\Sigma(\TT_2)$-coarse function of metric spaces then:
\begin{center}
$asdim_{\TT_1} X \le asdim_{(\TT_1, \TT_2)}(f) + asdim_{\TT_2}Y$.\\
\end{center}
\end{Thm}

Here is the global version:
\begin{Thm}\label{Hur}
If $f: X \to Y$ is a global $\Sigma(\TL_2)$-function of metric spaces then:
\begin{center}
$dim_{\TL_1} X \le dim_{(\TL_1, \TL_2)}(f) + dim_{\TL_2}Y$.\\
\end{center}
\end{Thm}

The following lemma appears in \cite{Hur} (Proposition $8.4.$) and it will be useful to apply the
 Hurewicz theorem to exact
sequences of finitely generated groups.

\begin{Lem}
If $1 \to K \to G \to H \to 1$ is an exact sequence and $G$ is a finitely generated group,
then there are word metrics $d_G$ on
$G$ and $D_H$ on $H$ such that $f: (G, d_G) \to (H, d_H)$ is $1$-Lipschitz and for any $m$-dimensional control function
$D_K$ on $K$ the function:\par
\begin{center}
$D_f(r_G, R_H):= D_K(r_G+2R_H) + 2R_H$\\
\end{center}
is an $m$-dimensional control function of $f$.
\end{Lem}

Then if $D_K \in \TT$ for some control semigroup $\TT$ we will have

\begin{Cor}
If $1 \to K \to G \to H \to 1$ is an exact sequence of groups so that $G$ is finitely generated then
\begin{center}
$asdim_{\TT} (G, d_G) \le asdim_{\TT}(K, d_G|_K) + asdim_{\TT}(H, d_H)$.\\
\end{center}
for any word metrics $d_G$ on $G$ and $d_H$ on $H$.
\end{Cor}

Finally, we will show the analogous version of the Hurewicz theorem for groups acting on spaces with finite
 semigroup-controlled asymptotic dimension.

Let $G$ be a group acting by isometries on a metric space $X$ and let $R>0$.
Given {\itshape $x_0 \in X$ the $R$-stabilizer of $x_0$
 }is defined by $W_R(x_0)= \{\gamma \in G: d(\gamma \dot x_0, x_0) \le R\}$.

The following theorem is the analogous to Theorem $8.9.$ in \cite{Hur}. The proof is completely similar.

\begin{Thm}
Let $\TT_1$ and $\TT_2$ be two control semigroups and
let $G$ be a finitely generated group acting by isometries on a metric space $X$ such that $asdim_{\TT_2} X < +\infty$.
Fix a point $x_0 \in X$. If there exists a $(\TT_1, \TT_2)$-function $D^{(\TT_1, \TT_2)}: \RR_+ \times \RR_+ \to \RR_+$
such that for all $R>0$ the function
$g: \RR_+ \to \RR_+$ defined by $g(r) = D(r, R)$ is a $(k, \TT_1)$-dimensional control function of $W_R(x_0)$, then:
\par
\begin{center}
$asdim_{\TT_1} G \le k + asdim_{\TT_2} X$.\\
\end{center}
\end{Thm}


\begin{thebibliography}{99}
\bibitem{Assouad}
P. Assouad, {\em Sur la distance de Nagata}, C. R. Acad. Sci. Paris S«er. I
Math. {\bf 294} (1982), no. 1, 31--34.

\bibitem{Bell Dran Beld}
G. ~Bell, A.~Dranishnikov, \emph{Asymptotic dimensión in Bedlewo}, preprint.

\bibitem{Brod-Dydak-Higes-Mitra}
N.~Brodskiy, J.~Dydak, J.~Higes, A.~Mitra, \emph{Nagata-Assouad
dimension via Lipschitz extensions}, preprint math.MG/0601226.

\bibitem{Brod Dydak Higes Mitra2}
N.~Brodskiy, J.~Dydak, J.~Higes, A.~Mitra, \emph{Dimension zero at all scales}, preprint.

\bibitem{Buyalo1}
S.~Buyalo, \emph{Asymptotic dimension of a hyperbolic space and
capacity dimension of its boundary at infinity}, Algebra i analis
(St. Petersburg Math. J.), {\bf 17} (2005), 70--95 (in Russian).

\bibitem{Dran Smith}
A.~Dranishnikov, J. Smith, \emph{On asymptotic Assouad-Nagata dimension}, preprint math.GR/0607143.

\bibitem{Dran-Zar}
A.~Dranishnikov, M.~Zarichnyi, \emph{Universal spaces for
asymptotic dimension}, Topology Appl. {\bf 140} (2004), no.2-3,
203--225.

\bibitem{Gro asym invar}
M. Gromov, \emph{Asymptotic invariants for infinite groups}, in
Geometric Group Theory, vol. 2, 1--295, G.Niblo and M.Roller,
eds., Cambridge University Press, 1993.

\bibitem{Hur}
N.~Brodskiy, J.~Dydak, M.~Levin, A.~Mitra, \emph{Hurewicz theorem for Assouad-Nagata dimension}, preprint math.GT/0605416.

\bibitem{Isbell book}
J.R.~Isbell, \emph{Uniform spaces}, Mathematical Surveys, No. 12
American Mathematical Society, Providence, R.I. 1964 xi+175 pp.

\bibitem{Lang-Sch Nagata dim}
U.~Lang, T.~Schlichenmaier, \emph{Nagata dimension, quasisymmetric
embeddings, and Lipschitz extensions}, IMRN International
Mathematics Research Notices (2005), no.58, 3625--3655.

\bibitem{Protasov}
I.~Protasov, \emph{Survey of Balleans}, preprint.

\bibitem{Roe lectures}
J. Roe, \emph{Lectures on coarse geometry}, University Lecture
Series, 31. American Mathematical Society, Providence, RI, 2003.


\end{thebibliography}
\end{document}